\documentclass[matzei, numbook, envcountsame]{svjour}
\usepackage{amssymb, amsmath, amsfonts}

\def\CM{Cohen-Macaulay }
\def\Spec{\mathop{\mathrm{Spec}}}

\def\height{\mathop{\mathrm{height}}}
\def\gr{{\mathop{\mathrm{gr}}}_m}
\def\Cl{\mathop{\mathrm{Cl}}}

\def\sym{\mathcal{R}_{\it{s}}}
\def\H{\mathop{\mathrm{Hom}_A}}
\def\HB{\mathop{\mathrm{Hom}_B}}

\begin{document}

\title{Multi-symbolic Rees algebras and strong F-regularity}

\author{Anurag K. Singh}

\institute{Department of Mathematics, University of Utah, 155 S. 1400 E., Salt
Lake City, UT 84112--0090, USA \\ 
\email{singh@math.utah.edu}}

\date{Received: August 10, 1998}

\maketitle

\begin{abstract}
Let $I$ be a divisorial ideal of a strongly F-regular ring $A$. K.-i.~Watanabe
raised the question whether the symbolic Rees algebra $\sym(I) = \oplus_{n \ge
0}I^{(n)}$ is Cohen-Macaulay whenever it is Noetherian. We develop the notion
of multi-symbolic Rees algebras and use this to show that $\sym(I)$ is indeed
\CM whenever a certain auxiliary ring is finitely generated over $A$.
\end{abstract}

\section{Introduction}

In \cite{Waantic} K.-i.~Watanabe raised the issue whether for a divisorial
ideal $I$ of a strongly F-regular ring $A$, the symbolic Rees algebra
$\sym(I)=\bigoplus_{n \ge 0}I^{(n)}U^n$ is Cohen-Macaulay whenever it is
Noetherian. Watanabe showed that this is true when $I$ is an {\it
anti-canonical ideal}\/ i.e., an ideal of pure height one which represents the
inverse of the class of the canonical module of $A$ in the divisor class group
$\Cl(A)$. In this paper we work in the more general setting of multi-symbolic
Rees algebras and as a corollary of our main result, Theorem \ref{main}, we
obtain the following positive answer to Watanabe's question:

\begin{theorem} 
Let $(A,m)$ be a strongly F-regular ring with canonical ideal $\omega$. Given
an ideal $I$ of $A$ of pure height one, choose $J$ of pure height one such that
$[I]+[J]+[\omega]=0$ in $\Cl(A)$. If the multi-symbolic Rees algebra
$\sym(I,J)$ is finitely generated over $A$, then $\sym(I)$ is Cohen-Macaulay.
\end{theorem}

The hypothesis that $A$ is strongly F-regular is indeed used in an essential
way: Watanabe has constructed an example of an F-rational ring $A$ with a
divisorial ideal $I$ such that the symbolic Rees algebra $\sym(I)$ is not 
Cohen-Macaulay, see \cite[Example 4.4]{Waantic}.

In general, of course, the symbolic Rees algebra $\sym(I)$ for a divisorial
ideal $I$ of a normal ring $A$ need not be Noetherian, e.g., if $A$ is the
coordinate ring of an elliptic curve, and $I$ is a prime ideal of height one,
which has infinite order in the divisor class group, $\Cl(A)$. However if we
specialize to the case when $A$ is F-rational, a two-dimensional example is
easily ruled out since, by a result of J.~Lipman, \cite{Li}, the divisor class
group of a two dimensional rational singularity (and hence by \cite{Smratsing}
of a two dimensional F-rational ring) is a torsion group. In  dimension three
the hypothesis that $A$ has rational singularities is no longer sufficient:
S.~D.~Cutkosky has shown that a symbolic Rees algebra over a three dimensional
ring with rational singularities need not be Noetherian, see \cite[Theorem
6]{dale}. It should be noted that if $A$ is a Gorenstein ring of dimension
three over $\mathbb C$ with rational singularities, then symbolic Rees algebras
at divisorial ideals of $A$ are finitely generated by \cite[Theorem
6.1]{kawamata}.

\section{Preliminaries}

Throughout our discussion all rings are commutative and have a unit element.
Unless stated otherwise, we shall assume our rings contain a field $K$ of
characteristic $p > 0$. We use the letter $e$ to denote a variable nonnegative
integer, and $q$ to denote the $e\,$th power of $p$. We denote by $F$ the
Frobenius endomorphism of $A$, i.e., $F(a) = a^p$. For a reduced ring $A$ of
characteristic $p > 0$, $A^{1/q}$ shall denote the ring obtained by adjoining
all $q\,$th roots of elements of $A$. The ring $A$ is said to be {\it
F-finite}\/ if $A^{1/p}$ is module-finite over $A$. Note that a finitely
generated algebra $A$ over a field $K$ is F-finite if and only if $K^{1/p}$ is
a finite field extension of $K$.

In the notation $(A,m)$, the ring $A$ is either a Noetherian local ring with
maximal ideal $m$, or an $\mathbb N$-graded ring with homogeneous maximal ideal
$m = \oplus_{i > 0}A_i$ which is finitely generated over a field $A_0 = K$. 

By a normal domain, we shall mean a Noetherian domain which is integrally
closed in its field of fractions.

Our references for the theory of tight closure are \cite{HHstrong},
\cite{HHjams}, \cite{HHbasec}, and \cite{HHjalg}. We next recall some
definitions and well known facts.

\begin{definition}
A ring $A$ is said to be F-pure if for all $A$-modules $M$, the Frobenius 
homomorphism $F: M \to F(M)$ is injective.
 
An F-finite domain $A$ is strongly F-regular if for every nonzero element $c
\in A$, there exists $q=p^e$ such that the $A$-linear inclusion $A \to A^{1/q}$
sending $1$ to $c^{1/q}$ splits as a map of $A$-modules.   
\end{definition}

Regular rings are strongly F-regular, and strongly F-regular rings are
Cohen-Macaulay. If $B$ is a strongly F-regular ring and $A$ is a subring which
is a direct summand of $B$ as an $A$-module, then $A$ is also strongly
F-regular.

Let $(A,m)$ be an F-finite local domain, and let $E=E_A(A/m)$ denote the
injective hull of the residue field $A/m$. K.~E.~Smith has shown that $A$ is
strongly F-regular if and only if the zero submodule of $E$ is tightly closed,
see \cite[Proposition 7.1.2]{Smthesis}. Hence if $\zeta \in E$ is a socle
generator, $A$ is strongly F-regular if and only if for every nonzero element 
$c \in R$, there exists a positive integer $e$ such that $cF^e(\zeta) \neq 0$.

\begin{definition}
Let $A$ be a normal domain and $I$ an ideal of pure height one. Then $I^{(n)}$ 
denotes the $n\,$th symbolic power of the ideal $I$, i.e., the reflexive hull 
of $I^n$. If $\mathcal{F}$ is the set of minimal primes of $I$, we have
$$
I^{(n)} = (\bigcap_{P \in \mathcal{F}}I^nA_P) \cap A.
$$
\end{definition}

\section{Multi-symbolic Rees algebras}

Let $(A,m)$ be a Noetherian normal local domain, and $I$ an ideal of pure
height one. The symbolic Rees algebra 
$$
\sym(I)=\bigoplus_{n \ge 0}I^{(n)}U^n \subseteq A[U]
$$ 
is an object that has been studied extensively. We generalize this construction 
to a finite family of ideals $I_1, I_2, \dots, I_k$, each of pure height one,
by defining
$$
\sym(I_1, I_2, \dots, I_k) = \bigoplus_{n_1, \dots, n_k \in \mathbb N} 
(I_1^{n_1} I_2^{n_2} \dotsm I_k^{n_k})^{**} 
\ U_1^{n_1}U_2^{n_2}\dotsm U_k^{n_k} 
$$ 
as a subring of the polynomial ring $A[U_1, \dots, U_k]$. Here $^*$ denotes 
the dual $\H(-, A)$. In this notation $(I_1^{n_1} \dotsm I_k^{n_k})^{**}$ 
is the reflexive hull of $I_1^{n_1} \dotsm I_k^{n_k}$.

\begin{proposition} 
Let $(A,m)$ be a normal domain and $I_1, \dots, I_k$ be ideals of $A$ of pure 
height one. Then the multi-symbolic Rees algebra $B =\sym(I_1, \dots, I_k)$ is 
a Krull domain. Hence the ring $B$ is a normal domain whenever it is Noetherian.
\end{proposition} 

\begin{proof}
Let $\mathcal{F}$ be the set of all minimal prime ideals of $I_1, \dots, I_k$. 
We then have
$$
\sym(I_1, \dots, I_k) 
 = (\bigcap_{P \in \mathcal{F}} A[I_1 A_P U_1, \dots, I_k A_P U_k]) 
 \cap A[U_1, \dots, U_k]
$$
which, being a finite intersection of Krull domains, is a Krull domain. 
\qed \end{proof}

\begin{proposition}
Let $(A,m)$ be a normal domain, and $I_1, \dots, I_k$ be ideals of $A$ of pure 
height one such that the multi-symbolic Rees algebra 
$B = \sym(I_1, I_2, \dots, I_{k-1})$ is Noetherian. Let 
$\tilde I_k=(I_kB)^{**}$ denote the reflexivization of $I_kB$ as a $B$-module, 
i.e., $^*$ denotes $\HB(-, B)$. Then there is a natural isomorphism
$$
\sym(\tilde I_k)= B \oplus \tilde I_k \oplus \tilde I_k^{(2)} \oplus \cdots
\cong \sym(I_1, I_2, \dots, I_k).
$$ 
\label{mainprop}
\end{proposition}

\begin{proof}
There is a natural inclusion
$$
B \oplus I_k U_k B \oplus I_k^2 U_k^2 B \oplus \cdots \to \sym(I_1,\dots, I_k).
$$
To obtain the isomorphism asserted, we need to verify that a reflexive 
$B$-module is reflexive when considered as an $A$-module. For this it 
suffices to verify that $B$ is a reflexive $A$-module, but this follows since 
$B$ is a direct sum of reflexive $A$-modules. 
\qed \end{proof}

This gives us the immediate corollary:

\begin{corollary}
Let $(A,m)$ be a normal domain, and $I_1, \dots, I_k$ ideals of  $A$ of pure
height one such that $B = \sym(I_1, I_2, \dots, I_k)$ is Noetherian.  Then $B$
arises by a successive construction of symbolic Rees algebras starting  with
the ring $A$.  
\label{iterate} \end{corollary}

\begin{theorem}
Let $(A,m)$ be a normal domain, and $I_1, \dots, I_k$ be ideals of  $A$ of pure
height one such that the multi-symbolic Rees algebra $B = \sym(I_1, I_2, \dots,
I_k)$ is Noetherian. Then the inclusion $A \subseteq B$ satisfies Samuel's PDE
condition, i.e., for a height one prime $P \in \Spec B$, we have $\height(P
\cap A) \le 1$. This gives a natural map of divisor class groups, $i: \Cl(A)
\to \Cl(B)$, which is an isomorphism. 

Furthermore if $A$ and $B$ are homomorphic images of regular local rings and
$\omega_A$ and $\omega_B$ denote the canonical modules of $A$ and $B$
respectively, we have the relation
$$
[\omega_B] = i([\omega_A] + [I_1] + \dots + [I_k]).
$$
\label{class}
\end{theorem}

\begin{proof}
The corresponding statements for symbolic Rees algebras (i.e., the case $k=1$) 
are covered by \cite[Lemma 4.3 (2), Proposition 4.4, Theorem 4.5]{GHNV}, see
also \cite[Proposition 2.6]{simis-trung}. The assertions here follow by 
combining these results with Corollary \ref{iterate} above. 
\qed \end{proof}

Determining when symbolic Rees algebras are Cohen-Macaulay is a subtle issue:
as remarked earlier, Watanabe has constructed examples where $A$ is a ring with
rational singularities and $I$ is an anti-canonical ideal, but the symbolic
Rees ring $\sym(I)$ is not  Cohen-Macaulay. However if the divisorial ideals
$I_1, \dots, I_k$ have finite order as elements of the divisor class group
$\Cl(A)$, we have the following extension of \cite[Theorem 4.1]{GHNV}:

\begin{theorem}
Let $(A,m)$ be a normal domain, and $I_1, \dots, I_k$ be ideals of pure height
one which have finite order as elements of the divisor class group $\Cl(A)$.
Then the multi-symbolic Rees algebra $B = \sym(I_1, \dots, I_k)$ is \CM if
and only if for all $n_i \in \mathbb N$ the ideals $(I_1^{n_1} \dotsm
I_k^{n_k})^{**}$ are maximal \CM $A$-modules.
\end{theorem}

\begin{proof}
Let $a_i$ denote the order of $[I_i]$ in $\Cl(A)$, and fix elements $x_i$ such 
that $I_i^{(a_i)} = x_i A$, for $1 \le i \le k$. The elements 
$x_1 U^{a_1}, \dots, x_k U^{a_k}$ form part of a system of parameters for $B$,
and it is easily verified that this is a regular sequence on $B$. Next note
that 
$$
B / (x_1 U^{a_1}, \dots, x_k U^{a_k}) = \bigoplus_{0 \le n_{ij} < a_i} 
(I_1^{n_{1j}} \dotsm I_k^{n_{kj}})^{**} U_1^{n_{1j}} \dotsm U_k^{n_{kj}} 
$$
and so $B$ is a \CM ring if and only every system of parameters for $A$ is a 
regular sequence on the ideals $(I_1^{n_1} I_2^{n_2} \dotsm I_k^{n_k})^{**}$
for all $n_1, \dots, n_k \in \mathbb N$.
\qed \end{proof}

\section{Examples}

\begin{example}
Consider the subring $A=K[ax, ay, bx, by]$ of the polynomial ring $K[a,b,x,y]$ 
and the height one prime ideals $P_i=(ax, ay)$ and $Q_j=(ax, bx)$ where 
$1 \le i \le n$ and $1 \le j \le m$. Then the multi-symbolic Rees algebra
$$
B = A(P_1, \dots, P_n, Q_1, \dots, Q_m)
$$ 
is isomorphic to the Segre product of two polynomial rings, 
$$
K[X_1, \dots, X_{n+2}] \#K[Y_1, \dots, Y_{m+2}].
$$
We have $\Cl(A) = \Cl(B) = \mathbb Z$, and 
$$
[\omega_B] = i( n[P] + m[Q]) = i ( (n-m)[P] ).
$$
In particular, $B$ is Gorenstein if and only if $n=m$.
\end{example}

\begin{example}
Let $A = K[X_1, \dots, X_m]^{(n)}$ denote the $n\,$th Veronese subring of the 
polynomial ring $K[X_1, \dots, X_m]$. We compute all multi-symbolic Rees
algebras over the ring $A$. The divisor class group of $A$ is $\Cl(A) =
{\mathbb Z}/n{\mathbb Z}$ and we fix as a generator, the height one prime ideal 
$$
P = (X_1K[X_1, \dots, X_m]) \cap A
$$
Divisorial ideals of $A$ are of the form $P^{(i)}$ up to isomorphism, and the
multi-symbolic Rees algebra $A( P^{(\alpha_1)}, \dots, P^{(\alpha_k)} )$ is
determined by the $k$-tuple of integers $\alpha_1, \dots, \alpha_k$. We claim
that $A( P^{(\alpha_1)}, \dots, P^{(\alpha_k)} )$ is isomorphic to the $n\,$th
Veronese subring of the polynomial ring
$$
K[X_1, \dots, X_m, X_1^{\alpha_1}U_1, \dots, X_1^{\alpha_k}U_k]
$$
where the variables $U_1, \dots, U_k$ have weight zero. To see this, note that 
by definition we have 
$$
A( P^{(\alpha_1)}, \dots, P^{(\alpha_k)} ) = \bigoplus_{n_i \ge 0} 
P^{ ( n_1 \alpha_1 + \dots + n_k \alpha_k ) } U_1^{n_1} \dotsm U_k^{n_k} 
$$
and that a monomial in $X_1, \dots, X_m$ is an element of
$P^{ ( n_1 \alpha_1 + \dots + n_k \alpha_k ) }$ precisely if it is a multiple 
of $X_1^{n_1 \alpha_1 + \dots + n_k \alpha_k}$ whose degree is a multiple of 
$n$.
\end{example}

\section{An application to tight closure theory}

Let $(A,m)$ be a strongly F-regular domain with canonical module $\omega$. In
\cite{Waantic} Watanabe showed that the anti-canonical symbolic Rees algebra
$\sym(I)$ is Cohen-Macaulay (in fact, strongly F-regular) whenever it is
Noetherian, and raised the question whether this is true for an arbitrary ideal
$I$ of pure height one. As an application of the construction of
multi-symbolic Rees algebras, we show that $\sym(I)$ is strongly F-regular,
and in particular is Cohen-Macaulay, whenever a certain auxiliary algebra is
finitely generated over $A$. Our main theorem is:

\begin{theorem}
Let $(A,m)$ be an F-finite normal ring with canonical ideal $\omega$. Given an
ideal $I$ of $A$ of pure height one, choose $J$ of pure height one such that
$[I]+[J]+[\omega]=0$ in the divisor class group $\Cl(A)$. Assume that the 
multi-symbolic Rees algebra $\sym(I,J)$ is finitely generated over $A$. If $A$ 
is F-pure, then the rings  $\sym(I)$ and $\sym(I,J)$ are also F-pure. If $A$ is
strongly F-regular,  then $\sym(I)$ and $\sym(I,J)$ are strongly F-regular, and
in particular, are Cohen-Macaulay.
\label{main}
\end{theorem}

\begin{proof}
Let $B=\sym(I) = \oplus_{i \ge 0}I^{(i)}U^i$ and $\tilde J = (JB)^{**}$ where
$^*$ denotes $\HB(-, B)$. If $R=\sym(I,J) \subseteq A[U,V]$, by 
Proposition \ref{mainprop} we have 
$$
R=\sym(\tilde J) = B \oplus \tilde J \oplus \tilde J^{(2)}  \oplus \cdots
$$ 
Setting $d = \dim A$, we have $\dim R = d+2$. Consider the maximal ideal of $B$,
$$
\mathfrak{m} = m+ IU +I^{(2)}U^2 + \cdots
$$ 
and the maximal ideal of $R$,
$$
\mathfrak{M} = \mathfrak{m} + \tilde JV + \tilde J^{(2)}V^2 + \cdots .
$$
In \cite[Theorem 2.2]{Waantic} Watanabe has computed the highest local 
cohomology module of a symbolic Rees ring, and furthermore determined the 
Frobenius action on it. Using this we have 
$$
H^{d+2}_{\mathfrak{M}}(R) \cong \bigoplus_{j<0}
                     H^{d+1}_{\mathfrak{m}}(\tilde J^{(j)})V^j.
$$
Again using Watanabe's result we get 
$$
H^{d+2}_{\mathfrak{M}}(R) \cong \bigoplus_{i<0, j<0}
          H^{d}_m( I^{(i)} J^{(j)}) U^i V^j.
$$
By Theorem \ref{class} and the fact that $[I]+[J]+[\omega]=0$ in $\Cl(A)$ we
see that $R = A(I,J)$ is quasi-Gorenstein, i.e., has a trivial canonical
module. Hence $H^{d+2}_{\mathfrak{M}}(R)$ is the injective hull of $R/
\mathfrak{M}$, and so the strong F-regularity or F-purity of $R$ can be
determined by studying the action of the Frobenius on 
$H^{d+2}_{\mathfrak{M}}(R)$. 

Let $\zeta$ be a socle generator of $H^d_m(\omega)$. Then the socle of
$H^{d+2}_{\mathfrak{M}}(R)$ is generated by $\zeta U^{-1}V^{-1}$. If $A$ is
F-pure, then $F(\zeta) \neq 0$ and so  
$$
F(\zeta U^{-1}V^{-1} ) = F(\zeta)U^{-p}V^{-p} \neq 0,
$$
by which $R$ is F-pure. Consequently $B = \sym(I)$, being a direct summand of 
$R$, is also F-pure.

Next assume that $A$ is strongly F-regular. Let $cU^nV^m \in R$ be a nonzero
element where $c \in (I^nJ^m)^{**}$. To show that $R$ is also strongly 
F-regular, it suffices to show that 
$(cU^nV^m) \ F^e(\zeta U^{-1}V^{-1} ) \neq 0$ for some positive integer $e$. 
Choosing a suitable multiple, if necessary, we may assume that $n=m$. We may
choose a canonical module $\omega$ for $A$ such that $(IJ)^{**} =
\omega^{(-1)}$. Then $c \in (I^nJ^n)^{**} = \omega^{(-n)}$. Since $A$ is 
strongly F-regular we
may choose $e$ such that $q = p^e > n$ and $c F^e(\zeta) \neq 0$ in 
$H^d_m(\omega^{(q)})$. But then $c F^e(\zeta) \neq 0$ as an element of
$H^d_m(\omega^{(q-n)})$. Hence
$$
cU^nV^n \cdot F^e( \zeta U^{-1}V^{-1} ) = c F^e(\zeta) (UV)^{n-q} 
       \in H^d_m(\omega^{(q-n)})(UV)^{n-q} 
$$
is nonzero, and so $R$ is strongly F-regular. Hence its direct summand $B$ is
also strongly F-regular, and therefore is Cohen-Macaulay.
\qed 
\end{proof}

\section{Rees rings}

So far in our discussion, we had been considering symbolic Rees rings at ideals
of pure height one. In this section we switch to the other extreme and consider
the Rees ring at the homogeneous maximal ideal $m$ of an $\mathbb N$-graded
normal ring $(A,m)$. (In this case the Rees ring $R=A[mT]$ agrees with the
symbolic Rees ring $\sym(m)$.) Although the relation between the properties of
$A$ and those of $R=A[mT]$ seems to be very mysterious, there is one case where
easy answers are available:

\begin{proposition} 
Let $(A,m)$ be an $\mathbb N$-graded normal ring which is generated by its 
degree one elements over the field $A_0=K$. Consider the Rees ring $R=A[mT]$.
If $A$ is strongly F-regular, F-pure, or a ring of characteristic zero with 
rational singularities, then the same is true for $R$.  
\end{proposition}

\begin{proof} 
Note that $R=A[mT]$ is isomorphic to the Segre product $A\# B$ where $B=K[S,T]$
is a polynomial ring in two variables. Consequently $R$ is a direct summand of
$A[S,T]$. If $A$ has rational singularities, then so does $A[S,T]$, and
consequently $R$ has rational singularities (in characteristic $0$) by 
Boutot's result, \cite{boutot}. Similarly if $A$ is strongly F-regular or
F-pure, the same is true for $A[S,T]$, and its direct summand, $R$. 
\qed  \end{proof}

In the following example $A$ is a normal monomial ring, i.e., a normal subring
of a polynomial ring which is generated by monomials. Consequently $A$ is
strongly F-regular but we shall see that the Rees ring $A[mT]$ fails to be
normal.

\begin{example} Consider the monomial ring 
$$ 
A = K[W^3X, \ X^3Y, \ Y^3Z, \ Z^3W, \ W^2X^2Y^2Z^2] \subseteq K[W,X,Y,Z] 
$$ 
where $K$ is a field. It is not difficult to see that $A$ is isomorphic to the 
the hypersurface 
$$
K[U_0, U_1, U_2, U_3, U_4]/(U_0^2 - U_1 U_2 U_3 U_4)
$$ 
and is a normal ring. Consequently by the main result of \cite{Ho}, $A$ is a
direct summand of a regular ring, and so is strongly F-regular. Take the Rees
ring be $R=A[mT]$. The element $W^2X^2Y^2Z^2T^2$ is in the fraction field of
$R$, although it is not in $R$ itself. However 
$$
(W^2X^2Y^2Z^2T^2)^2 = (W^3XT) (X^3YT) (Y^3ZT) (Z^3WT) \in R
$$
and so $R$ is not normal. Furthermore, when the characteristic of the field $K$ 
is $2$, it is easily verified that $R$ is not F-pure, although the ring $A$ 
is F-pure.
\end{example}

In the next example $A$ is an F-rational hypersurface, but the Rees ring
$A[mT]$, while being Gorenstein and normal, is not F-rational.

\begin{example}  Let $A = K[W,X,Y,Z]/(W^2+X^3+Y^6+Z^7)$. Then $A$ is
F-rational whenever the characteristic of K is $p \ge 7$. We show that the
Rees ring $R = A[mT]$, while being Gorenstein and normal, is not F-rational.

First note that the Rees ring $R$ is Gorenstein. This holds, for example, 
by \cite[Theorem 1.2]{goto-shimoda} since the associated graded ring 
$$
\gr(R) \cong K[W,X,Y,Z]/(W^2)
$$ 
is Gorenstein with $a$-invariant $a(\gr(R)) = -2$. 

We next examine $R$ on the punctured spectrum. For $f \in m$, the localization
$R_f \cong A_f[T]$ is a polynomial ring over $A_f$. For an element $fT$ with 
$f \in m$, note that 
$$
R_{fT} \cong K[\frac{w}{f}, \ \frac{x}{f}, \ \frac{y}{f}, \ \frac{z}{f}, \ f, \
fT, \ \frac{1}{fT}].
$$
Examining these localizations as $f$ ranges through the set $\{w,x,y,z\}$, we
can see that $R$ is indeed normal. 

To see that $R$ is not F-rational, take $f = z$ above and let 
$U_1=\frac{w}{z}$, $U_2=\frac{x}{z}$ and $U_3=\frac{y}{z}$. Then 
$R_{zT} \cong S[zT, 1/zT]$ where
$$
S = K[U_1, \ U_2, \ U_3, \ Z]/(U_1^2 + U_2^3Z + U_3^6Z^4 + Z^5).
$$
It suffices to show that $S$ is not F-rational. Consider the grading on $S$
where the variables $U_1, \ U_2, \ U_3, \ Z$ have weights $15, 8, 1, 6$
respectively. The $a$-invariant of $S$ is easily computed to be $a(s) = 0$, 
and so $S$ cannot be F-rational.
\end{example}

\begin{acknowledgement}
It is a pleasure to thank Mel Hochster for several valuable discussions.
\end{acknowledgement}

\end{document}